\documentclass[preprint,12pt]{elsarticle} 

\usepackage[utf8]{inputenc}
\usepackage{amssymb}
\usepackage{amsmath}
\usepackage{listings}
\usepackage{xcolor}

\begin{document}

\begin{frontmatter}


\title{A study of concurrent multi-frontal solvers for modern massively parallel architectures}

\author[label1]{Jan Trynda}
\author[label1]{Maciej Wo\'{z}niak}
\author[label2]{Sergio Rojas}
\address[label1]{AGH University of Science and Technology \\ Institute of Computer Science \\ al. A Mickiewicza 30, 30-059 Krak\'{o}w, Poland \\ email: macwozni@agh.edu.pl } 
\address[label2]{Instituto de Matem\'aticas, Pontificia Universidad Cat\'olica de Valpara\'iso. Valpara\'iso, Chile\\ email: sergio.rojas.h@pucv.cl}

\begin{abstract}


Leveraging Trace Theory, we investigate the efficient parallelization of direct solvers for large linear equation systems.
Our focus lies on a multi-frontal algorithm, and we present a methodology for achieving near-optimal scheduling on modern massively parallel machines.
By employing trace theory with Diekert Graphs and Foata Normal Form, we rigorously validate the effectiveness of our proposed solution.
To establish a strong link between the mesh and elimination tree of the multi-frontal solver, we conduct extensive testing on matrices derived from the Finite Element Method (FEM).
Furthermore, we assess the performance of computations on both GPU and CPU platforms, employing practical implementation strategies.
\end{abstract}

\begin{keyword}
Trace theory \sep Direct solver \sep Concurrent computations \sep Matrix factorization
\end{keyword}

\end{frontmatter}


\lstset
{ 
    basicstyle=\fontsize{9}{11}\selectfont\ttfamily,
    numbers=left,
    stepnumber=1,
    showstringspaces=false,
    tabsize=1,
    breaklines=true,
    breakatwhitespace=false,
}

\tableofcontents

\section{Introduction}
\label{sec:introduction}


In modern simulations utilizing methods like the Finite Difference Method (FDM), Finite Element Method (FEM), and others, a two-phase approach is commonly followed.
First, the matrices representing the discrete system are obtained, and then the resulting linear or nonlinear algebraic equation system is solved.
Solving the systems of linear equations can be accomplished through multifrontal direct solvers such as MUMPS \cite{MUMPS1,MUMPS2,MUMPS3,MUMPS4}, SuperLU \cite{SUPERLU1,SUPERLU2}, or PaStiX \cite{PASTIX}, which are known as industry standards.
Alternatively, iterative solvers like PCG \cite{PCG}, GMRES \cite{banas}, or the Additive Schwarz method \cite{domain_decomposition} can be employed.

Iterative solvers offer significant advantages, particularly their high scalability on large clusters, enabling the efficient handling of extensive systems of equations.
Compared to direct solvers, iterative solvers typically require fewer computational resources regarding time and memory.
However, they do have drawbacks.
One major issue is the potential for severe convergence problems, especially when the matrix's conditioning is large and/or efficient preconditioners are unavailable.
Additionally, in certain scenarios, iterative solvers may perform worse than direct solvers.
This is particularly evident when solving multiple matrices that share the same set of rows and columns, as in mesh-based methods when local grid refinements are performed \cite{paszynska2015quasi,ICCSGALOIS,aboueisha2014dynamic}.

It is important to note that direct solvers play a vital role within multiple iterative solvers and are essential in various applications.
Different solvers are required for specific applications such as elasticity \cite{el2010iterative}, electromagnetism \cite{hiptmair1998multigrid}, and fluid dynamics \cite{arnold2000multigrid}, air pollution propagation \cite{OLIVER201347, podsiadlo2021parallel}, as well as different numerical methods.
Despite the convergence problems associated with iterative solvers, direct solvers remain valuable in many scenarios.
Furthermore, iterative solvers may be slower than direct solvers in certain situations, highlighting the need to carefully consider the choice of solver based on the specific characteristics of the problem at hand.

\subsection{State-of-the-art algorithm for direct solving sparse linear systems of equations: Multi-frontal solver}


The multi-frontal solver is a widely recognized algorithm for directly solving sparse linear systems of equations.
Duff and Reid initially proposed it in their sequential version of the algorithm \cite{duff1984multifrontal,duff1983multifrontal}, which utilizes the nested dissection algorithm \cite{khaira1992nested}.
The core idea behind the multi-frontal solver is to recursively solve the problem by traversing the nodes of the elimination tree.

The elimination tree serves as the input for the algorithm, where the frontal matrices are located in the tree's leaves.
The solver algorithm follows the elimination tree from the leaves to the root, solving subproblems at each node.
This process is repeated until reaching the root of the elimination tree.
Once the root problem is solved, a recursive backward substitution is performed, descending from the root to the leaves, to obtain the final solution.

Among the various direct methods available for solving linear equations, such as LU factorization, QR factorization, and singular value decomposition \cite{van1996matrix}, the multi-frontal solver commonly employs LU factorization, also known as Gaussian elimination.
LU factorization decomposes the original matrix A into the product of a lower triangular matrix L and an upper triangular matrix U.

The fundamental principle of the LU factorization algorithm is to transform the original system of equations into an equivalent system where the coefficient matrix is factored into L and U.
This factorization allows for efficient forward and backward substitutions to obtain the solution to the system.
The multi-frontal solver incorporates LU factorization as a key step in its solution process, making it a powerful and effective algorithm for directly solving sparse linear systems.

\subsection{Architecture}
\label{sec:architecture}

\begin{figure}[ht]
    \centering
    \includegraphics[width=1.0\textwidth]{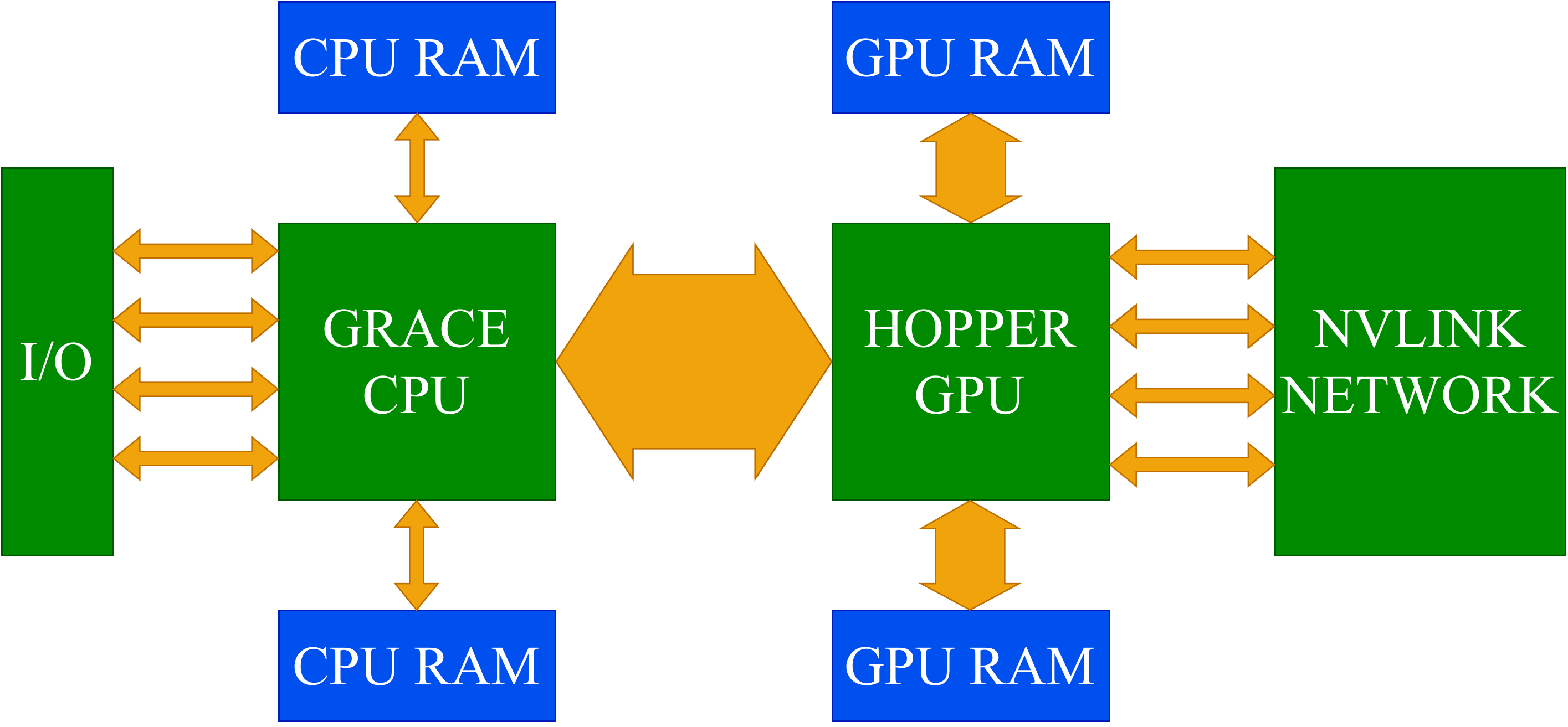}
    \caption{Simplified of a single Nvidia Grace-Hopper Superchip compute unit.
    Current state-of-the-art for future supercomputers and data centers.}
    \label{fig:cluster}
\end{figure}

State-of-the-art supercomputers, exemplified in Figure \ref{fig:cluster}, are sophisticated multi-level hierarchical hybrid systems \cite{cyfronet,stampede,Summit}.
They consist of classical nodes (servers) interconnected via high-speed networks, such as Infiniband, and utilize the Message Passing Interface (MPI) logic for communication.
Each server can house multiple multi-core CPUs sharing with a shared RAM, while also accommodating massively parallel co-processors like GPUs.
These GPUs possess dedicated memory with a hierarchical structure \cite{CUDAmemory}, distinct from the CPU's memory.
To optimize hardware utilization, promote eco-friendly computing, and reduce carbon emissions, it is essential to develop concurrent algorithms tailored specifically for these systems, ensuring efficient resource utilization and a smaller environmental impact.

\subsection{Trace theory}
\label{sec:tracetheory}


Various methods can be employed to verify concurrent computations formally.
One widely utilized approach is \textit{Trace Theory} \cite{TraceTheory}, which encompasses the Foata Normal Form (FNF) \cite{BookOfTraces} and Diekert dependency graphs.
These techniques are particularly valuable for enabling the parallel implementation of algorithms on highly parallel architectures, such as GPU clusters.
Moreover, Trace Theory is a foundation for near-optimal scheduling and theoretically verifies parallel algorithm correctness.
Formal verification also leverages other methods, including \textit{Petri Net} \cite{PetriNet}, \textit{Process Calculi} \cite{ProcessCalculi}, \textit{TLA$^+$} \cite{lamport2002specifying}, and the \textit{Actor Model} \cite{ActorModel}.
These approaches contribute to the formal verification landscape, offering diverse tools and techniques for ensuring the correctness and reliability of concurrent computations.

\subsection{Schedule of research}
\label{sec:scheduleofresearch}


This article focuses on multi-frontal parallel direct solver algorithms designed specifically for modern hybrid memory clusters that feature multiple GPUs per node.
Our primary emphasis is on the LU factorization method.
While there are alternative algorithms available, such as cuSOLVER or rocSOLVER \cite{cusolver,rocsolver,GHYSELS2022102897}, it is worth noting that they lack a formal proof of optimality and do not support cluster configurations.

The trace-theory-based analysis of concurrency presented in this article is not limited to the factorization algorithm alone but can be extended to other contemporary direct factorization methods applicable to general sparse matrix structures.
In the interest of clarity and focusing on the trace theory methodology, we have chosen to demonstrate its application in the context of classical FEM, associating the elimination tree with the computational mesh.
This simplification allows for a more focused presentation while still showcasing the potential of trace theory in a specific domain.
Our approach divides the solver algorithm into multiple indivisible tasks, considering fundamental mathematical operations like summation, multiplication, and more.
To establish the execution order of these tasks, we conducted a thorough analysis of their dependencies using trace theory.
By introducing dependency relations, we identified classes of tasks that can be executed concurrently.
These task classes conform to the Foata Normal Form (FNF), which we achieved through trace theory.

To validate our methodology, we extensively tested the implementation of the solver on an NVIDIA CUDA architecture.
This rigorous testing on CUDA ensures the reliability and effectiveness of our approach, demonstrating its suitability for practical applications.

While several articles have explored similar topics, it is essential to note the distinctions.
Some of these works have applied Trace Theory to multi-frontal solver algorithms \cite{kuznik2012graph,obrok2010graph,ICCSGALOIS}.
However, they treat operations within a single front as atomic without further decomposing into basic operations such as multiplication or summation.
Additionally, these articles often introduce algorithms closely tied to specific graph grammars and FEM computations.

Other works focus on specific types of parallelism, such as shared memory or distributed memory approaches \cite{wozniak2015computational,paszynski2010parallela,wozniak2014computational,LOS201799}.
While they contribute valuable insights, they do not address formal concurrency verification specifically for modern massively parallel machines.

It is worth highlighting that, as of now, none of the published works or algorithms formally verify concurrency for modern massively parallel machines.
This gap underscores the novelty and significance of our research in providing a formal verification approach that addresses the challenges of concurrency on these advanced computational platforms.

\subsection{Structure of the article}
\label{sec:strucutre}


The article is organized as follows.
Section~\ref{sec:formulation} introduces the model problem that serves as the basis for the subsequent benchmarking process.
Section~\ref{sec:concurency_model} delves into the multi-frontal algorithm and explores how Trace theory is leveraged to develop a concurrent algorithm variant.
Section~\ref{sec:numres} presents the results of several numerical experiments conducted to evaluate the algorithm's performance.
A detailed analysis of the findings is also provided.
Lastly, Section~\ref{sec:conclusions} offers a comprehensive summary of the paper, highlighting the key contributions and potential future directions for research.

\section{Formulation of the model problems}
\label{sec:formulation}

\subsection{The discrete $L^2$-projection}
\label{sec:L2}


We present the $L^2$-projection into a finite element space as a motivational example to illustrate the proposed methodology in this work.
However, it is important to note that the methodology can be applied to any finite element problem that involves the inversion of an invertible symmetric matrix.

Let $\Omega \subset \mathbb{R}^d$ be an open and bounded polyhedron.
Using standard notation for Hilbert spaces, the discrete $L^2$-projection problem can be stated as follows: Given a function $f \in L^2(\Omega)$ and a finite-dimensional space $V_h \in L^2(\Omega)$, we aim to find a discrete solution $u_h \in V_h$ that satisfies the equation:
\begin{equation}\label{eq:L2_proj}
(u_h,v_h) = (f,v_h), \quad \forall v_h \in V_h,
\end{equation}
where $(\cdot,\cdot)$ represents the $L^2$ inner product.
In essence, we are seeking the best approximation of $f$ within the space $V_h$ in the $L^2$ sense.

Consider a basis ${\phi_i}_{i=1}^n$ for the finite element space $V_h$, and let $\mathbf{u}=(u_1,\dots,u_n)^T$ be the vector of coefficients of $u_h$ with respect to this basis.
Thus, the solution $u_h$ to problem~\eqref{eq:L2_proj} can be expressed as:
\begin{equation}
u_h=\sum_{i=1}^n u_i , \phi_i.
\end{equation}
It is well-known that solving problem~\eqref{eq:L2_proj} leads to a linear system of equations in the form:
\begin{equation}\label{eq:Linear_system}
M \mathbf{u} = \mathbf{b},
\end{equation}
where $M$ is the symmetric and positive definite matrix, referred to as the mass matrix, with entries $M_{ij} = (\phi_j, \phi_i)$ for $i,j = 1, \dots, n$, and $\mathbf{b}$ is the load vector with entries $b_i=(f,\phi_i)$ for $i=1,\dots,n$.

\subsection{Schur complement}
\label{sec:schur}


The Schur complement method is a remarkable and widely employed technique for solving linear systems.
It provides an efficient approach to handling such systems by leveraging the inherent properties and structure of the matrices involved.
By decomposing the original matrix into submatrices and isolating the Schur complement, which is itself a well-defined matrix, we can simplify the problem and optimize the computational process.

Consider a linear system represented as $B\mathbf{u}=\mathbf{b}$, where $B$ is a matrix that can be decomposed into submatrices.
Specifically, decomposed into the form:
\begin{equation}\label{eq:Bmat}
\begin{aligned}
B =
\begin{bmatrix}
C & D \\ E & -F
\end{bmatrix},
\end{aligned}
\end{equation}
where $C$ and $F$ are square matrices of dimensions $q$ and $r$, respectively.
If $C$ is invertible, then the Schur complement method consists of performing a partial LU factorization on the submatrix $C$, resulting in the following factorization:
\begin{equation}
\begin{aligned}
B =
\begin{bmatrix}
I & 0 \\ EC^{-1} & I
\end{bmatrix}
\begin{bmatrix}
C & 0\\ 0 & -(F+EC^{-1}D)
\end{bmatrix}
\begin{bmatrix}
I & C^{-1}D \\ 0 & I
\end{bmatrix}.
\end{aligned}
\end{equation}
The term $-(F+E^TC^{-1}D)$ is known as the Schur complement.

In practical applications, matrix $C$ often represents a well-conditioned and computationally efficient system component, while $F$ may pose conditioning or computational cost challenges. By applying the Schur complement method, we can solve the linear system efficiently by first solving for the matrix $C$ and then utilizing this solution to eliminate the dependence on $C$ in the Schur complement matrix. The Schur complement plays a vital role in resolving partial differential equations (PDEs) through residual minimization on dual discrete norms, which is a remarkable example of its application. In a seminal work by authors in \cite{CohDahWelM2AN2012}, it was proven that this type of discretization leads to the inversion of a matrix with a specific structure, as given in equation \eqref{eq:Bmat}, where $C$ is a symmetric and positive definite matrix, $E = D^T$, and $F=0$. Notably, the Schur complement, in this case, is also a symmetric and positive definite matrix. This key observation has paved the way for the development of various methods, including the DPG method \cite{DemGopCMAME2010, DemGopNMPDE2011, demkowicz2012class, zitelli2011class} (with a more recent overview available in \cite{DemGopBOOK-CH2014}), the AS-FEM \cite{calo2020adaptive, cier2021automatically, rojas2021goal, kyburg2022incompressible} (along with an isogeometric version discussed in \cite{los2021dgirm}), the iGRM \cite{CALO2021113214, LOS2020213, LOS2021200}, among others. These methods have leveraged the properties of the Schur complement to design efficient and robust techniques for solving PDEs, opening up new avenues in computational modeling and simulation.

\subsection{Multi-frontal algorithm}
\label{sec:MF-alg}


The frontal solver algorithm, initially introduced by M. Irons in \cite{irons1970frontal}, revolutionized the field of linear systems of equations.
This algorithm, rooted in Gaussian elimination, was originally implemented in Fortran.
Over time, it gained significant recognition and became a widely adopted standard for industrial and research applications.

I.S. Duff and J.K. Reid pioneered one notable advancement to the frontal solver algorithm \cite{Duff:1983a, Duff:1983b}.
Their work introduced the concept of the multi-frontal solver, which proved to be a pivotal modification.
The multi-frontal solver was initially devised for indefinite symmetric matrices and later extended to unsymmetric matrices.
This modification enabled the simultaneous elimination of multiple fronts, giving rise to the name ``multi-frontal''.
The multi-frontal solver's ability to handle complex matrix structures and exploit parallelism has contributed to its widespread adoption and applicability in various domains.
Its impact on solving linear systems of equations has been substantial, making it a cornerstone in the field of numerical methods.

\begin{figure}[ht]
    \centering
    \includegraphics[width=1.0\textwidth]{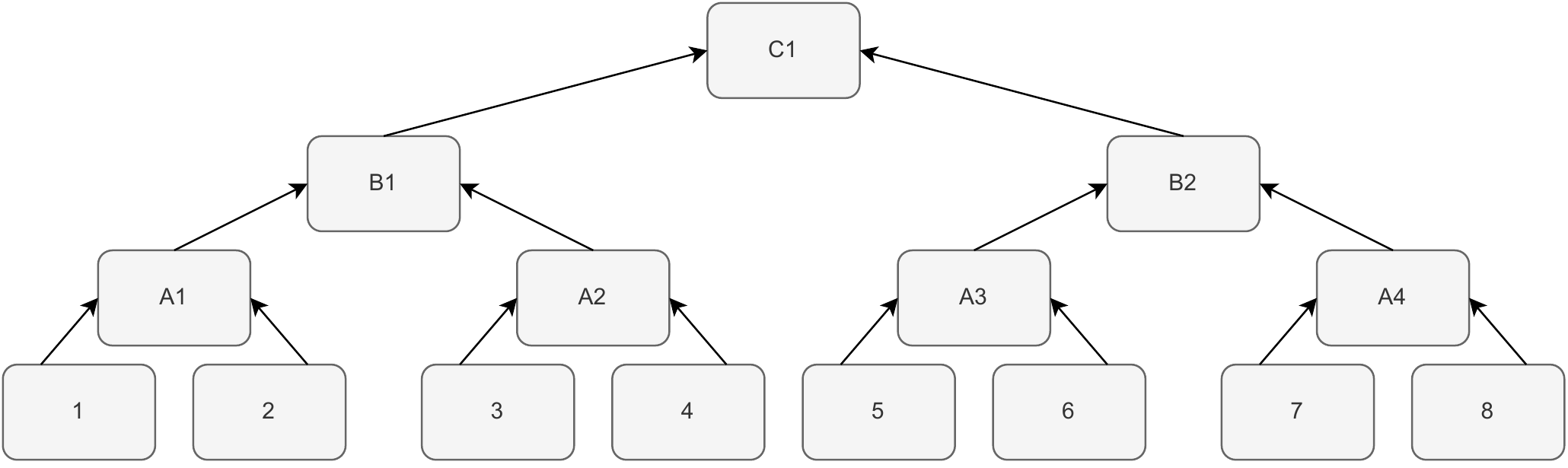}
    \caption{Simple example of a dependency tree for the (multi-)frontal method}
    \label{fig:mf_dep_tree}
\end{figure}

The multi-frontal method is a powerful approach for solving systems of linear equations, where $Ax=b$, with $A$ being a sparse matrix and $b$ a vector.
This method organizes the matrix factorization operations by dividing the matrix into smaller dense submatrices and performing partial factorization.

The dependencies between these partial factorizations, known as the Schur complement, can be represented as a tree structure (Figure \ref{fig:mf_dep_tree}).
To compute a node in the tree, all its child nodes must be computed first, ensuring proper dependencies are maintained.
This simultaneous computation of nodes based on their dependencies has been extensively discussed in \cite{paszynski2010parallel,geng1997parallel,Duff:1983a}.

Figure \ref{fig:mf_dep_tree} presents an assembly tree with eight nodes (fronts).
Notably, computations in fronts $A1$, $A2$, $A3$, and $A4$ can be performed simultaneously.
Furthermore, front $B1$ only requires the completion of fronts $A1$ and $A2$, while front $B2$ requires fronts $A3$ and $A4$.
This dependency structure is commonly referred to as the elimination tree.

As the partial elimination progresses, adjacent fronts are merged into larger ones.
This allows the processing of incompletely assembled parts of a front at higher levels of the tree.
The process is repeated at each level of the tree until reaching the root, where the largest, dense, fully assembled matrix is solved.

The algorithm for the multi-frontal method is presented in Table \ref{tab:mult_method_pseudocode}, summarizing the steps involved in the factorization process.

Overall, the multi-frontal method's effectiveness lies in exploiting parallelism and efficiently handling sparse matrices through a hierarchical factorization approach.

\begin{table}[ht]
\centering
\begin{tabular}{c}
\noindent\begin{minipage}{1.05\linewidth}
\begin{lstlisting}[escapeinside={(*}{*)}, frame = single, framexleftmargin=20pt]
function multi_frontal_method(matrix):  
    fronts = divide_into_fronts(matrix) 
    while(len(fronts) > 1): # for each elimination tree level
        for front in fronts: # for each tree node
            compute_schur_complement(front) 
        end for
        fronts = join_fronts(fronts) # based on the elimination tree
    end while
    compute_dense_matrix(fronts[0])
    return fronts[0]
end function
\end{lstlisting}
\end{minipage}
\end{tabular}
\caption{Multi-frontal method pseudocode}
\label{tab:mult_method_pseudocode}
\end{table}

\subsection{Elimination Tree}
\label{sec:elim_tree}


In our case, applying the multi-frontal solver has been specifically tailored to the Finite Element Method (FEM) problem.
Extensive research, such as the work by Paszyńska et al.~\cite{paszynska2015quasi}, has been dedicated to connecting the elimination tree with the computational domain in FEM.
To showcase the simplicity of our approach, we generate fronts as elemental matrices.

In the context of our methodology, the assembly of the elimination tree follows a bottom-up approach.
The algorithm can be divided into four distinct steps.
\begin{enumerate}
    \item Generate elemental matrices (fronts) by partitioning the FEM mesh into distinct elements.
    \item Assemble the left-hand-side sparse matrix by combining the elemental matrices based on mesh connectivity.
    \item Sort the fronts in a specific order that reflects their origin within the mesh.
    \item Pair the fronts together to form a binary tree structure, ensuring that the remaining front remains unchanged in the case of an odd number of fronts.
    \item Repeat the previous step with the newly created or merged fronts until only the root problem remains.
\end{enumerate}

\begin{figure}[ht]
    \centering
    \includegraphics[width=0.9\textwidth]{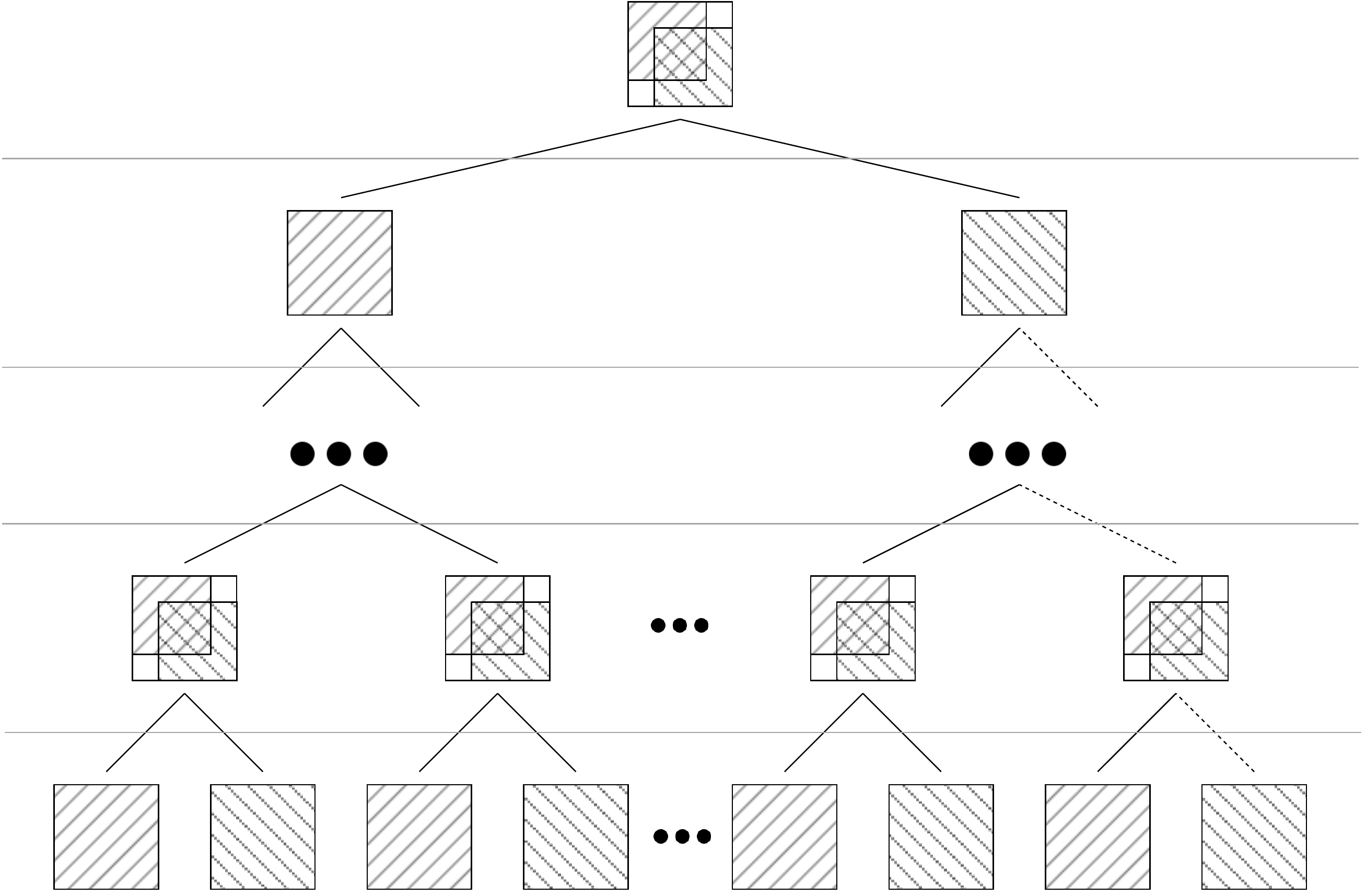}
    \caption{Visualisation of elimination tree.}
    \label{fig:elim_tree}
\end{figure}

The assembled elimination tree forms a balanced binary structure, with the potential exception of the farthest right element on each level (excluding the root).
Figure \ref{fig:elim_tree} visually represents the elimination tree, where each square represents a front.

The assembly process involves merging smaller fronts from lower levels to form larger fronts.
However, there can be instances where this rule does not apply.
For instance, a dashed line in the far right section of the tree highlights such an exception.
The farthest right element at each level may be formed from two or just one front from a lower level.
A simple identity relation exists between the two fronts in the latter case.

By isolating and examining the two-to-one relations, we can focus on the core process of merging smaller fronts to construct the elimination tree.

The composition of the newly assembled front can be described as follows:
\begin{itemize}
\item One part that is shared between both of the original fronts.
Our attention will be directed toward the yet-to-be-assembled portion of this front during this step of the elimination tree.
\item Two parts that are exclusive to only one of the original fronts.
\item Two parts that are derived from any previous front, consisting solely of zeros.
\end{itemize}

\section{Concurency model}
\label{sec:concurency_model}


Formally verifying concurrent computations involves creating a concurrency model using various methods.
It requires sophisticated concurrency analysis and meticulous construction of parallel algorithms.
To construct the Diekert-Mazurkiewicz trace model, we define an alphabet of tasks and establish a symmetric dependency relation between them.
In our model, each task from the alphabet occurs precisely once, ensuring a clear representation of dependencies.
To simplify the dependency relation and Diekert Graph, we introduce an additional type of task acting as a barrier.
While not directly originating from the simplest algorithm, this task is crucial in streamlining the concurrency analysis.
By employing this methodology, we can comprehensively discuss the full concurrency of the multi-frontal solver algorithm.
By carefully defining task alphabets, capturing dependencies, and incorporating barrier tasks, we enable a thorough exploration of parallel execution.

The alphabet of tasks for the multi-frontal factorization algorithm consists of the following four tasks.
Inside the definition of the tasks, we use ${z}$ as the global equivalent of index ${p}$, ${y}$ as the global equivalent of index ${j}$, and ${x}$ as the global equivalent of row index ${i}$.

\begin{enumerate}
    \item ${T^{1;f}_{\delta;\gamma^1}}$ - Division task: Divides a value ${b^{f}_{\delta}}$, where ${\delta=(p,j)}$, and $j$ is the element column index, on the selected pivot row $p$ by the current pivot ${b^{f}_{\delta^1}} $, where ${\delta^1=(p,p)}$ itself for the current front $f$.
    It computes the result of ${t^{1,f}_{\delta;\gamma^1} = b^{f}_{\delta} / b^{f}_{\delta^1}}$.
    This task is performed once for each element to the right (${j > p}$) of the pivot in the same row and the results are used for subsequent multiplication tasks.
    The ${\gamma^1}$ index points out the location of the elements in the global matrix.
    \item ${T^{2;f}_{\beta;\alpha}}$ - Multiplication task: Multiplies the result ${t^{1,f}_{\delta;\gamma^1}}$ of the corresponding division task ${T^{1;f}_{\delta;\gamma^1}}$, where $\delta = (p,j)$ by an element $b^{f}_{\delta^2}$, where ${\delta^2 = (i,p)}$ and ${\beta=(p,i,j)}$, where ${i}$ and $j$ are the elements row indexes of the selected pivot column ${p}$ for the current pivot $b^{f}_{\delta^1} $, with $\delta^1=(p,p)$ and front $f$.
    It computes the result of ${t^{2,f}_{\beta;\alpha} = t^{1,f}_{\delta;\gamma^1} * b^{f}_{\delta^2} }$.
    This task is performed once for each element below (${i > p}$) the pivot element in the same column.
    Each result will be used exactly once for subtraction, covered in the next task type.
    The index ${\alpha = (z,x,y)}$ points out the element's location in the global matrix.
    \item ${T^{3;f}_{\beta;\alpha}}$ - Subtraction task: Subtracts the results ${t^{2,f}_{\beta;\alpha}}$ of the corresponding multiplication task ${T^{2;f}_{\beta;\alpha}}$, where ${\alpha = (z,x,y)}$, from an element ${a_{\gamma}}$, where ${\gamma = (x,y)}$ of the global matrix and assigns the newly calculated value to this element.
    It computes the result of ${a_{\gamma} = a_{\gamma} - t^{2,f}_{\beta;\alpha}}$.
    This task is performed once for each multiplication task result.
    \item ${T^{4}_{\gamma}}$ - Assertion task: Checks the condition that the element ${a_{\gamma}}$, where ${\gamma = (x,y)}$ is ready to be used in calculations by verifying if it has reached its final value.
    This task performs an assertion operation on ${a_{\gamma}}$.
    It serves to ensure that all assignments to the elements have been completed.
\end{enumerate}

To summarize, each task in the multi-frontal factorization algorithm is represented by one or two upper subscripts, divided into two groups by a semicolon.
The first group indicates the task type, while the second group (optional) represents the current front index, denoted as $f$.
Additionally, tasks may have one or two subscripts at the bottom.
It is crucial to note that the execution of a specific task depends on completing tasks that provide the required input.
Thus, the algorithm ensures that the necessary dependencies are satisfied before executing a particular task.

\subsection{Set of dependencies}

In this section, we define the alphabet of tasks $\Sigma$ and the dependencies between them denoted by $D$.
The defined alphabet will consist of the following: 
\begin{itemize}
    \item division tasks for each nonzero element in the pivot row (where ${j > p}$), for each pivot, for each front,
    \item multiplication tasks for each nonzero element in the pivot column (where ${i > p}$), for each division task result,
    \item subtraction tasks for each multiplication task result,
    \item assertion tasks for each global matrix element (not only those initially nonzero),
\end{itemize}  
Therefore, a final form of the alphabet will be defined as follows:
\begin{equation}
\Sigma = \{
    T^{1;f}_{\delta;\gamma^1},
    T^{2;f}_{\beta;\alpha},
    T^{3;f}_{\beta;\alpha},
    T^{4}_{\gamma}
    \}
\label{eq:alphabet}
\end{equation}
Here index $f$ belongs to the set of all fronts in elimination tree.
Indexes $z$, $y$, and $x$ belong to the sub-matrix of the front $f$.
Now we can describe subsets of dependencies to simplify the construction.

The division depends on two assertions concerning the elements it uses for calculations.
Therefore, we can represent those dependencies as follows:
\begin{equation}
J_1 = \{
    (T^{4}_{\gamma^1}, T^{1;f}_{\delta;\gamma^1}),
    (T^{4}_{\gamma^2}, T^{1;f}_{\delta;\gamma^1})
    \}
\end{equation}

Multiplication depends on two tasks.
The first task is the corresponding division task that produces the value used during the current calculations.
The second task is the assertion that concerns the second element used for the calculations.
Therefore, we can represent those dependencies as follows:
\begin{equation}
J_2 = \{
    (T^{1;f}_{\delta;\gamma^1}, T^{2;f}_{\beta;\alpha}),
    (T^{4}_{\gamma^3}, T^{2;f}_{\beta;\alpha})
    \}
\end{equation}

Subtraction depends on the corresponding multiplication and multiple subtractions that
assign values to the same element ${a_{\gamma}}$ in the global matrix.
Therefore, we can represent those dependencies as follows:
\begin{equation}
J_3 = \{
    (T^{2;f}_{\beta;\alpha}, T^{3;f}_{\beta;\alpha}),
    (T^{3;f'}_{\beta';\alpha'}, T^{3;f}_{\beta;\alpha})
    \}
\end{equation} 
where $\alpha'=(z', x, y)$, and $( z' < z ) $ or $( z' = z \textrm{ and } f' < f )$

The assertion depends on all the subtractions that assign values to the same element ${a_{\gamma}}$ in the global matrix.
Therefore, we can represent those dependencies as follows:
\begin{equation}
J_4 = \{
    (T^{3;f'}_{\beta';\alpha'}, T^{4}_{\gamma})
    \}
\end{equation}
where $\alpha' = (z', x, y)$

Finally, we define dependencies between tasks from the alphabet $\Sigma.$
\begin{equation}
D = J^+ \cup (J^+)^{-1} \cup I_{\Sigma}
\label{eq:dependency}
\end{equation}
where
\begin{equation}
J = J_1 \cup J_2 \cup J_3 \cup J_4
\end{equation}

The primitives described above define the monoid of the traces for the problems under consideration.
We denote $J$ as the edges in the Diekert dependency graph \cite{TraceTheory}, which will be illustrated later in Figures \ref{fig:diekert_main}, \ref{fig:diekert_xy}, and \ref{fig:diekert_zz}.

After building the primitives of the trace monoid, i.e., the alphabet of tasks (\ref{eq:alphabet}) and the dependency relation (\ref{eq:dependency}), we define the pseudocode allowing us to compute matrix factorization, presented in Table \ref{tab:solver_pseudocode}.
The dependencies recorded in this algorithm determine only the sequence of operations in one string representing the desired trace.
By utilizing the task alphabet $\Sigma$ (\ref{eq:alphabet}), the dependency relation (\ref{eq:dependency}), and executing the pseudocode outlined in Tables \ref{tab:solver_pseudocode} and \ref{tab:elim_pseudocode}, we can calculate the Diekert dependency graph.
This graph is helpful for the appropriate and efficient scheduling of tasks in a heterogeneous computing environment.

\subsection{Application of trace theory to multi-frontal solver}
\label{sec:trace-mfs}


This section describes the methodology employed to generate the Diekert Dependency Graph (DDG) and the Foata Normal Form (FNF) in the context of the concurrent factorization of sparse matrices using the multi-frontal solver.
The DDG comprehensively represents all computational tasks involved in the factorization process, capturing their dependencies.
By utilizing the DDG and FNF, we can identify distinct groups of tasks known as Foata classes, which facilitate the practical implementation of concurrent computations.

To simplify the trace and dependency relation construction, we focus on a case where the sparse matrix is directly associated with the Finite Element Method (FEM) mesh.
The multi-frontal solver algorithm, presented in pseudocode in Table \ref{tab:solver_pseudocode}, assumes input divided into fronts based on the local contributions of each mesh element to the sparse matrix.

At each level of the elimination tree, the algorithm iterates over each front, merging it towards the tree's root while computing the Schur complement.
It is important to note that certain parts of a front may belong to multiple other fronts simultaneously.
In such cases, we select pivots that belong to only one front at a given time.
We reorder the front accordingly, placing the pivots to be computed at the current level first, followed by pivots shared by multiple fronts, and finally, pivots computed in previous levels of the elimination tree.
This reorganization allows for efficient computations.

The algorithm performs a division task within each pivot's row element.
Additionally, a multiplication task and subsequent subtraction task are performed for each element in the pivot column, with the resulting values assigned to the global matrix.
Finally, the algorithm computes the remaining front at the root level.

By following this methodology, we can effectively manage the concurrency of tasks in the multi-frontal solver, achieving efficient and accurate factorization of sparse matrices in a concurrent computing environment.

\begin{table}[ht!]
\centering
\begin{tabular}{c}
\noindent\begin{minipage}{1.05\linewidth}
\begin{lstlisting}[escapeinside={(*}{*)}, frame = single, framexleftmargin=20pt]
function multi_frontal_solver(fem_input, b_vector):  
    fronts = get_fronts(fem_input) 

    # For each elimination tree level
    while(len(fronts) > 1):
        # For each front
        for the front in fronts:
            # Compute what can be computed on this front
            compute_front(front)
        end for
        # Proceed to higher elimination tree level
        fronts = join_adjacent_fronts(fronts)
    end while
    compute_front(fronts[0])
    # Final factorized matrix
    final_matrix = fronts[0] 

    # Apply the B vector to return the result 
    return get_solution(final_matrix, b_vector) 


function compute_front(front):
    front = reorder(front)
    p = 0
    # As long as the selected pivot is in this front only and has not been computed yet
    while front[p, p].front_membership_count == 1 and !front[p, p].computed:
        for j in range(j, len(front)):
            # Division task
            t_1 = front[p,j] / front[p,p] 
            for i in range(i, len(front)):
                # Multiplication task
                t_2 = t_1 * front[i, p]
                # Subtraction task
                front[i, j].assign(front[i, j] - t_2)
            end for
        end for
        p += 1
    end while
    return front
\end{lstlisting}
\end{minipage}
\end{tabular}
\caption{Multi-frontal solver pseudocode}
\label{tab:solver_pseudocode}
\end{table}

\begin{table}[ht!]
\centering
\begin{tabular}{c}
\noindent\begin{minipage}{1.05\linewidth}
\begin{lstlisting}[escapeinside={(*}{*)}, frame = single, framexleftmargin=20pt]
function generate_elimination_tree(global_matrix):
    fronts = divide_into_FEM_based_fronts(global_matrix) 
    fronts = sort_fronts_based_on_lowest_index(fronts) 
    while(len(fronts) > 1):
        fronts = join_adjacent_fronts(fronts)
    end while
    return fronts[0] 
end function

function join_adjacent_fronts(fronts):
    new_fronts = [] 
    for i in range(0, len(fronts - 1), step=2): 
        new_fronts.append(join_fronts(fronts[i], fronts[i+1]))
    end for
    if len(fronts)%2 == 1: 
        new_fronts.append(fronts[len(fronts) - 1])
    end if
    return new_fronts
end function
\end{lstlisting}
\end{minipage}
\end{tabular}
\caption{Elimination tree generation pseudocode}
\label{tab:elim_pseudocode}
\end{table}

\begin{figure}[ht!]
    \centering
    \includegraphics[width=0.7\textwidth]{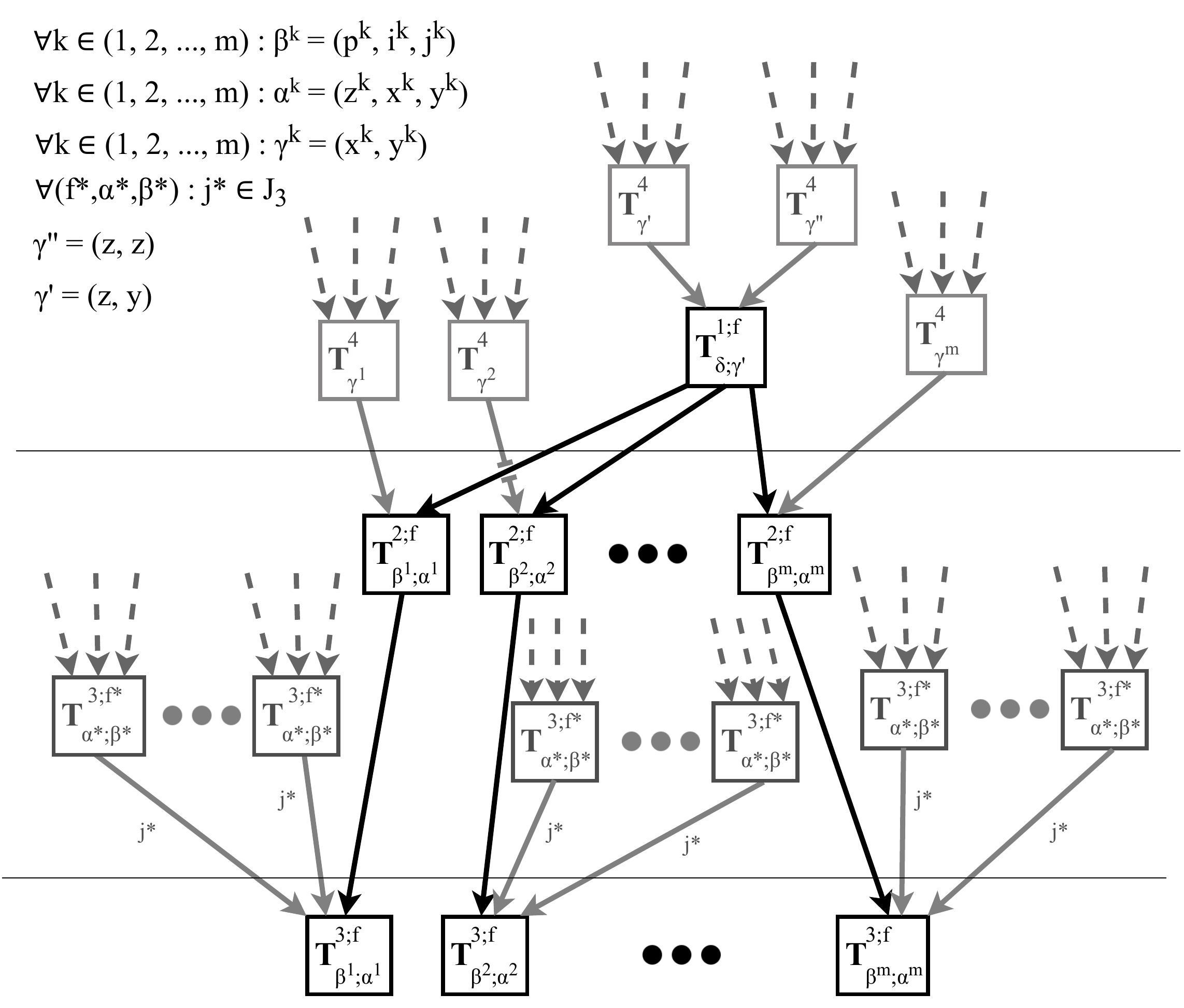}
    \caption{Main part of the Diekert Graph}
    \label{fig:diekert_main}
\end{figure}

We divided the Diekert graph into multiple recursive sub-graphs.
The first sub-graph focuses on operations within a single front, including division, multiplication, and subtraction tasks.
Figure \ref{fig:diekert_main} illustrates this sub-graph, with each task represented by a square with a signature.
Arrows depict the happen-before relation, denoted as ${J}$.

\begin{figure}[ht!]
    \centering
    \includegraphics[width=0.7\textwidth]{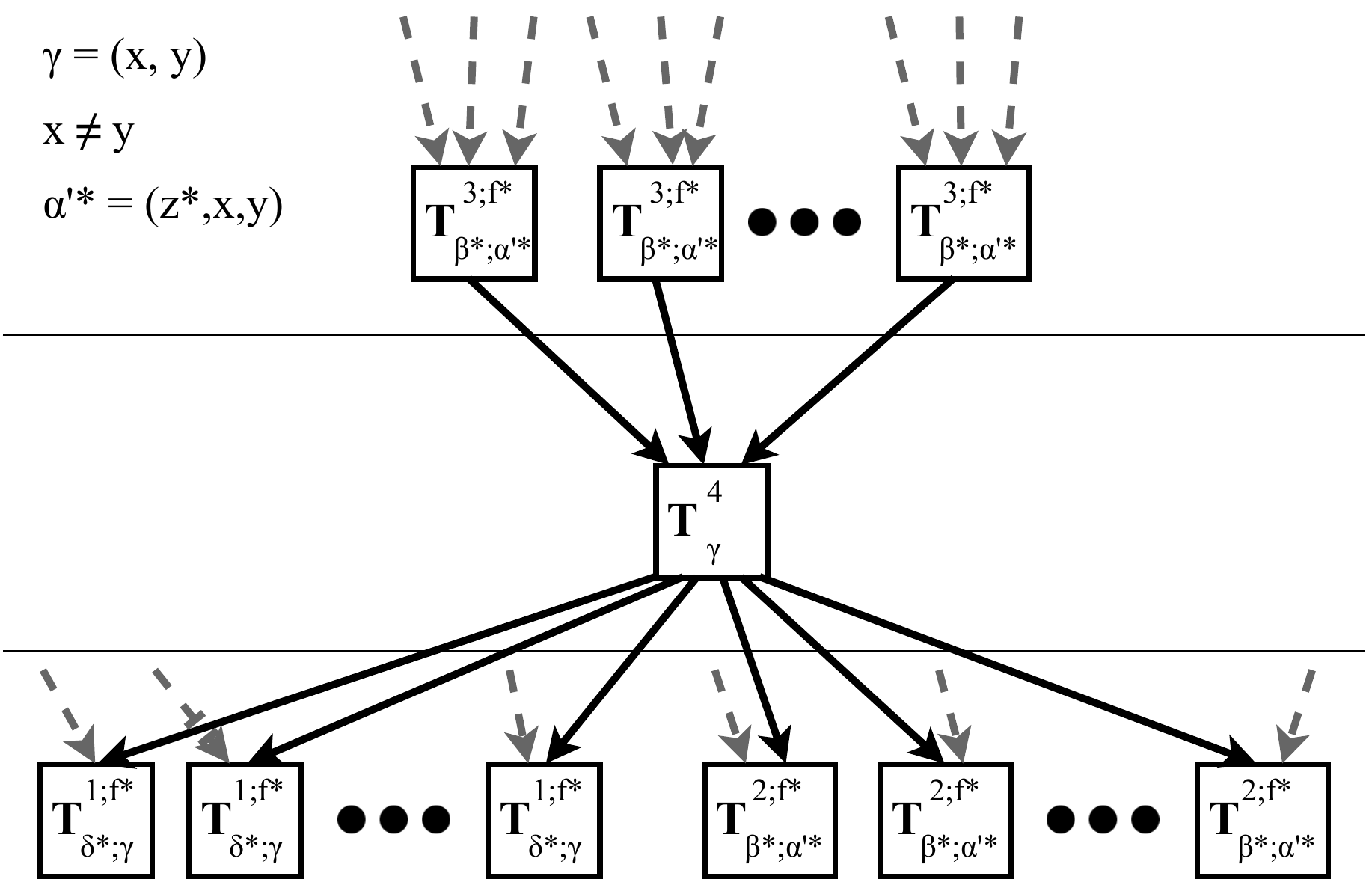}
    \caption{Diekert Graph part focused on a type 4 task regarding element off the diagonal}
    \label{fig:diekert_xy}
\end{figure}

\begin{figure}[ht!]
    \centering
    \includegraphics[width=0.7\textwidth]{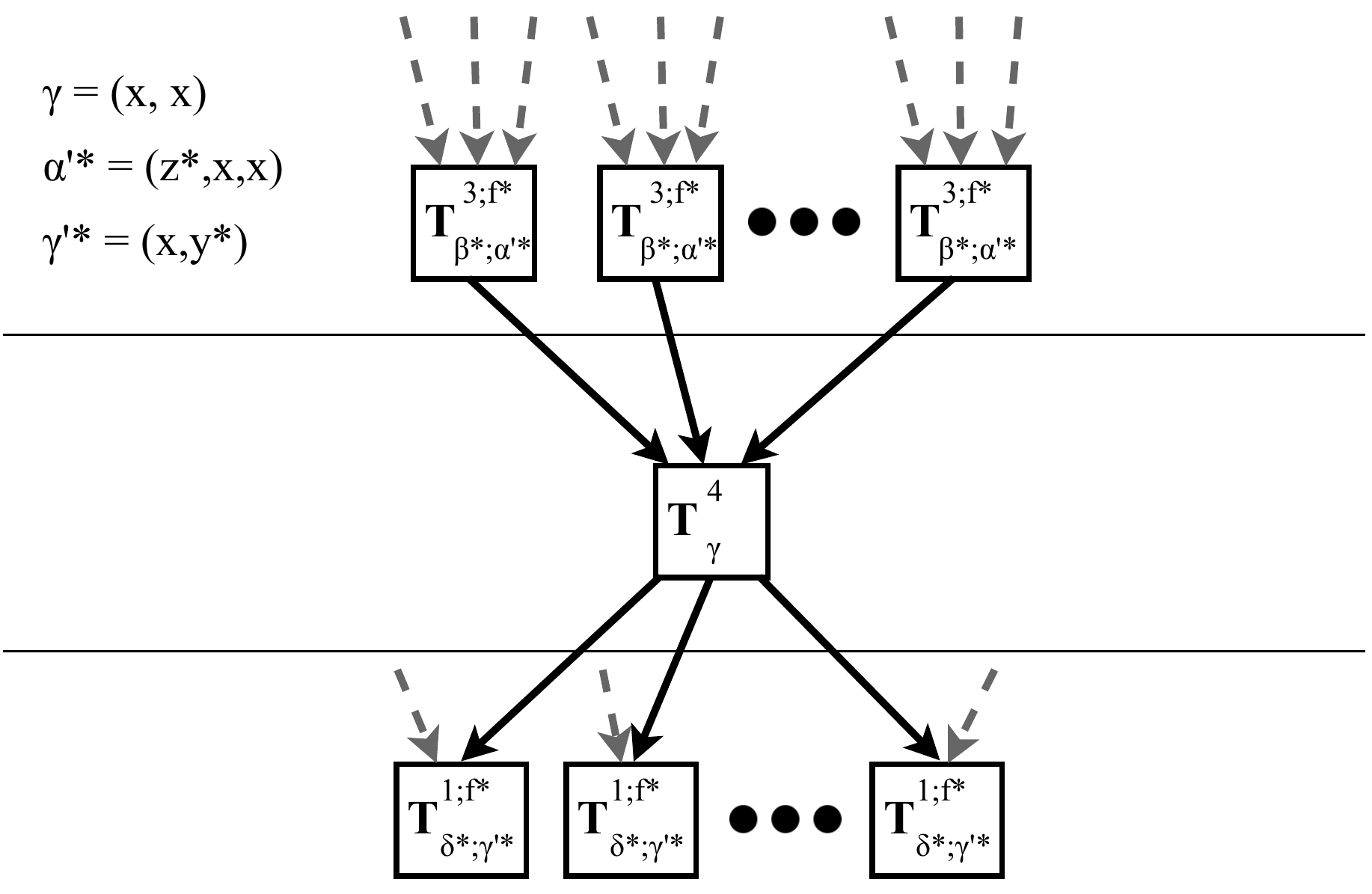}
    \caption{Diekert Graph part focused on a type 4 task regarding element on the diagonal}
    \label{fig:diekert_zz}
\end{figure}

Within this sub-graph, we can identify three types:
\begin{enumerate}
\item Sub-graph consisting of a single division task.
\item Sub-graph consisting of multiple multiplication tasks.
\item Sub-graph consisting of multiple subtraction tasks.
\end{enumerate}
The part of the graph involving assertion tasks requires two separate sub-graphs for different cases.
Figure \ref{fig:diekert_xy} represents the case with an off-diagonal element, while Figure \ref{fig:diekert_zz} represents the case with an element on the diagonal.
The tasks and relations in these sub-graphs are presented in the same manner as in the previous part of the graph.

Introducing these sub-graphs simplifies the Diekert graph by replacing many-to-many relations with many-to-one relations.
Consequently, we have a single relation in the first part of the Diekert graph.


\subsection{Scheduling algorithm}
\label{sec:fnf}


To achieve comparable scheduling performance to the classical algorithm on massively parallel shared-memory machines, we adopt the Foata-Normal-Form (FNF) approach \cite{traces}.
The Diekert dependency graphs (discussed in Section \ref{sec:trace-mfs}) illustrate the consecutive Foata classes for each scenario of the multi-frontal solver.
Within a particular Foata class, the tasks can be executed in any order, and the completion of the entire preceding Foata class is sufficient to initiate the computation of the next class.
This strategy ensures the absence of deadlocks, facilitates high-quality scheduling, and eliminates the need for intra-class synchronization.

To simplify the construction of the FNF, we first define groups of tasks that always belong to the same FNF class.
These task groups serve as a means to streamline the definition of FNF classes and enhance the overall organization of the algorithm.

For each front ($f$), each pivot, and each element in the pivot row ($\delta$, $\gamma^1$ - global index), we can define three groups of tasks:
\begin{equation}
    {p^{1;f}_{\delta;\gamma^1} = \{T^{1;f}_{\delta;\gamma^1}\}}
\end{equation}
\begin{equation}
    {p^{2;f}_{\delta;\gamma^1} = \{  T^{2;f}_{\beta^{1};\alpha^{1}},  T^{2;f}_{\beta^{2};\alpha^{2}}, \dots, T^{2;f}_{\beta^{m};\alpha^{m}} \}, \ \forall{l} \in (1, 2, \dots, m):  \beta^{l} = p, i^{l}, j } 
\end{equation}
\begin{equation}
    {p^{3;f}_{\delta;\gamma^1} = \{  T^{3;f}_{\beta^{1};\alpha^{1}},  T^{3;f}_{\beta^{2};\alpha^{2}}, \dots, T^{3;f}_{\beta^{m};\alpha^{m}} \}}
\end{equation}
We should also define an additional group of tasks for each global matrix element.
This is a synchronization/barrier group.
\begin{equation}
    {p^{4}_{\gamma} = \{T^{4}_{\gamma}\} }
\end{equation}
When defining the task groups, we can assume that each task belongs to exactly one group, corresponding to a single Foata class.
The Foata classes will be constructed by merging multiple task groups together.
For example, the tasks ${p^{1;f}_{\delta;\gamma^1}}$, ${p^{2;f}_{\delta;\gamma^1}}$, and ${p^{3;f}_{\delta;\gamma^1}}$ will be part of consecutive Foata classes in that specific order.

Based on Figures (\ref{fig:diekert_main}) through (\ref{fig:diekert_zz}), we can outline a general procedure for creating subsequent Foata classes, which involves the following tasks. 
\begin{itemize}
\item Class $q_1$, where $q \in \{1,5,9,\dots,4w+1\}$
\begin{equation}
\{p^{4}_{\gamma*}, p^{4}_{\gamma*}, \dots, p^{4}_{\gamma*}\}
\label{eq:class1}
\end{equation}

\item Class $q_2$, where $q \in \{2,6,10,\dots,4w+2\}$
\begin{equation}
\{p^{1;f^*}_{\delta^*;\gamma^*}, p^{1;f^*}_{\delta^*;\gamma^*}, \dots, p^{1;f^*}_{\delta^*;\gamma^*}\}
\label{eq:class2}
\end{equation}

\item Class $q_3$, where $q \in \{3,7,11,\dots,4w+3\}$
\begin{equation}
\{p^{2;f^*}_{\delta^*;\gamma^*}, p^{2;f^*}_{\delta^*;\gamma^*}, \dots, p^{2;f^*}_{\delta^*;\gamma^*}\}
\label{eq:class3}
\end{equation}

\item Class $q_4$, where $q \in \{4,8,12,\dots,4w+4\}$
\begin{equation}
\{p^{3;f^*}_{\delta^*;\gamma^*}, p^{3;f^*}_{\delta^*;\gamma^*}, \dots, p^{3;f^*}_{\delta^*;\gamma^*}\}
\label{eq:class4}
\end{equation}
\end{itemize}
Here, ${*}$ represents a placeholder for any (multi-)index value.
The notation used does not consistently refer to the same value, resulting in potential discrepancies. For instance, when considering an (multi-)index ${w}$, the value represented by ${w*}$ in one equation may not be equivalent to the value represented by ${w*}$ in other sections of the same or different equation.

We know that ${\forall{k}}$ if class $k$ takes the form
\begin{equation}
\{p^{1;f^*}_{\delta^*;\gamma^*}, p^{1;f^*}_{\delta^*;\gamma^*}, \dots, p^{1;f^*}_{\delta^*;\gamma^*}\}
\label{eq:class_k}
\end{equation}
then, the two following classes will have the form of
\begin{itemize}
\item Class $k+1$
\begin{equation}
\{p^{2;f^*}_{\delta^*;\gamma^*}, p^{2;f^*}_{\delta^*;\gamma^*}, \dots, p^{2;f^*}_{\delta^*;\gamma^*}\}
\label{eq:class_k1}
\end{equation}
\item Class $k+2$
\begin{equation}
\{p^{3;f^*}_{\delta^*;\gamma^*}, p^{3;f^*}_{\delta^*;\gamma^*}, \dots, p^{3;f^*}_{\delta^*;\gamma^*}\}
\label{eq:class_k2}
\end{equation}
\end{itemize}
where, if and only if, class ${k}$ contains element ${p^{1;f^*}_{\delta^*;\gamma^*}}$, then class ${k+1}$ contains element  ${p^{2;f^*}_{\delta^*;\gamma^*}}$ and class ${k+2}$ contains element ${p^{3;f^*}_{\delta^*;\gamma^*}}$.

The first type of Foata classes (\ref{eq:class1}) is responsible for assertions/synchronization.
Subsequent Foata classes follow this:
\begin{itemize}
\item The second type of Foata classes (\ref{eq:class2}) is responsible for divisions tasks.
\item The third type of Foata classes (\ref{eq:class3}) is responsible for multiplication based on results from the preceding second type Foata class (\ref{eq:class2}).
\item The fourth type of Foata classes (\ref{eq:class4}) is responsible for subtraction based on results from the preceding third type Foata class (\ref{eq:class3}).
\end{itemize}

All the mentioned tasks are executed on a uniform architecture, which allows us to anticipate similar execution times for each task within a specific Foata class.
This uniformity enables the efficient scheduling of tasks within a Foata class using a simple bag-of-tasks scheduling method.

We start the computation with Foata class 1 and ensure that all tasks within the class are completed before moving on to the next Foata class.
This sequential execution of Foata classes guarantees the correct order of computations and dependencies.

By adopting this simplified scheduling approach based on the Foata Normal Form (FNF), we achieve near-optimal performance in practical applications, even though no theoretical proof exists to support it.
This approach balances efficiency and simplicity, allowing for straightforward implementation and effective utilization of computational resources.

\section{Numerical results}
\label{sec:numres}


We conducted computational performance evaluations of the proposed algorithm on both CPU and GPU architectures.
Three modes of computation were considered: single-core CPU, multi-core CPU, and GPU.
The execution time was measured for the sequential factorization algorithm on the CPU and the concurrent integration algorithm on a shared memory CPU with four cores (8 threads).
The experiments were conducted on a Linux laptop equipped with an Intel Core i7-7700HQ processor and 24 GB RAM.
A Google Colab node with an Nvidia Tesla A100 GPU featuring 6912 CUDA cores was utilized for the GPU computations.

\begin{figure}[ht!]
    \centering
    \includegraphics[width=0.8\textwidth]{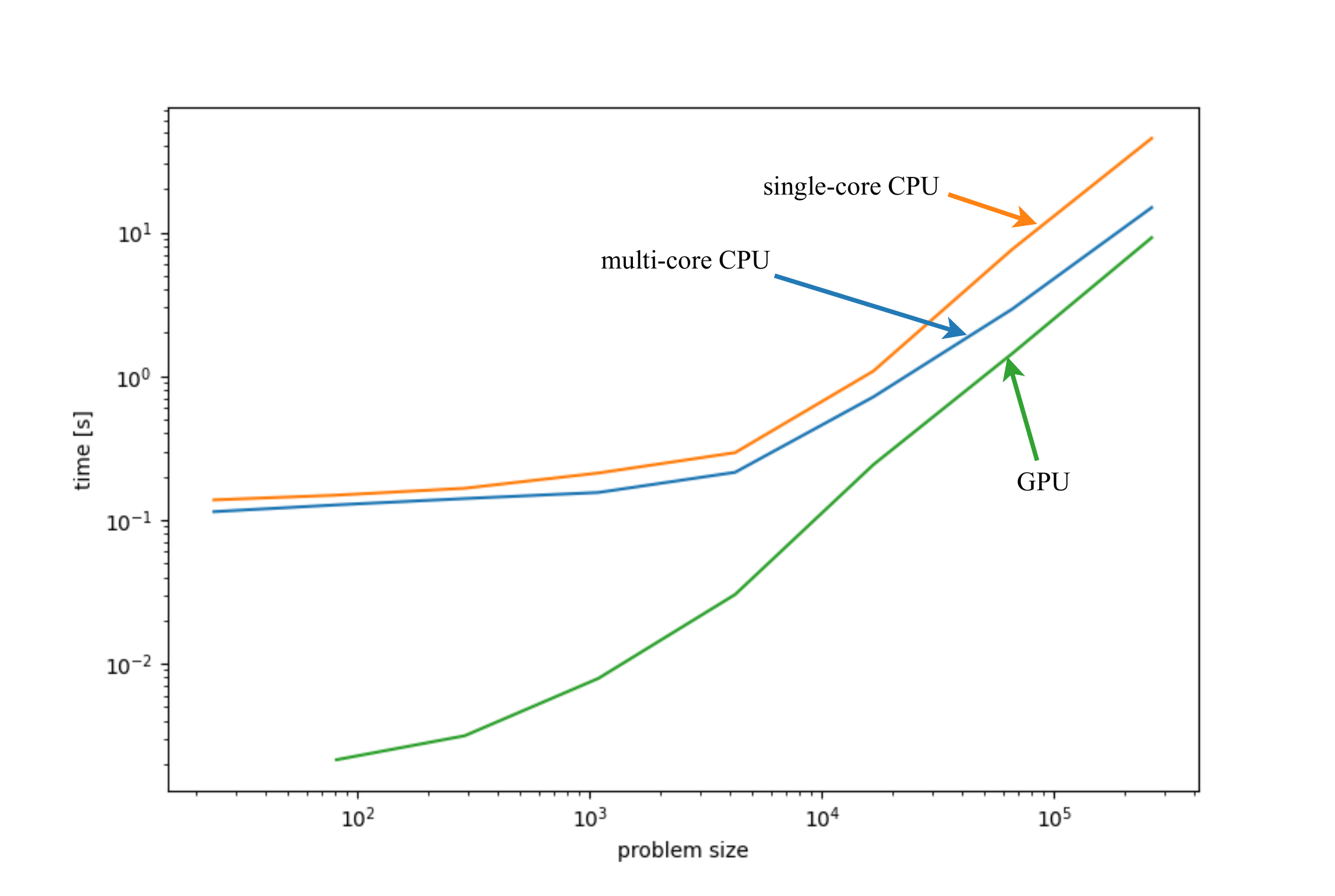}
    \caption{Comparison of CPU and GPU factorization time of sparse matrices for linear basis functions.}
    \label{fig:res_01}
\end{figure}

\begin{figure}[ht!]
    \centering
    \includegraphics[width=0.8\textwidth]{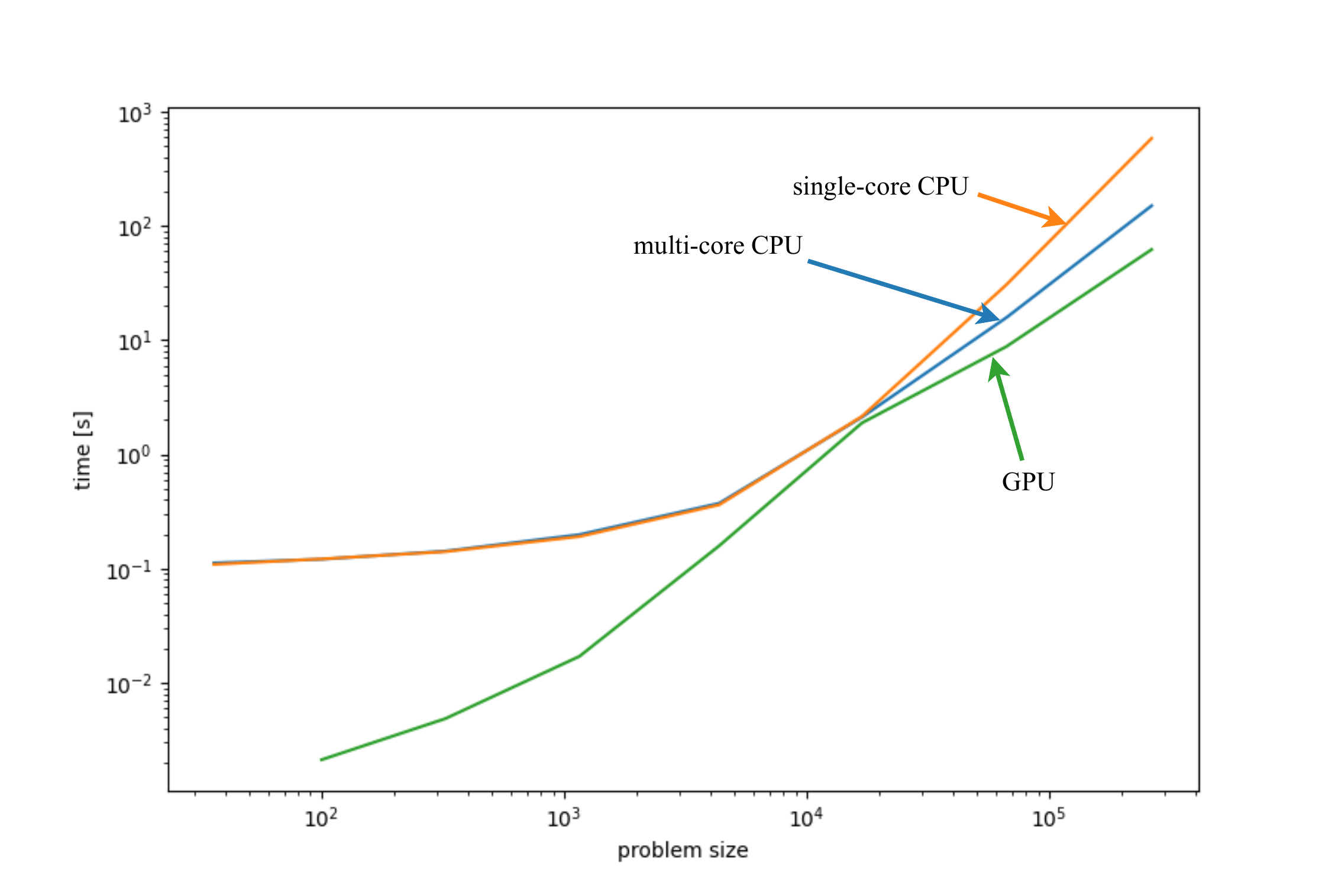}
    \caption{Comparison of CPU and GPU factorization time of sparse matrices for quadratic basis functions.}
    \label{fig:res_02}
\end{figure}

\begin{figure}[ht!]
    \centering
    \includegraphics[width=0.8\textwidth]{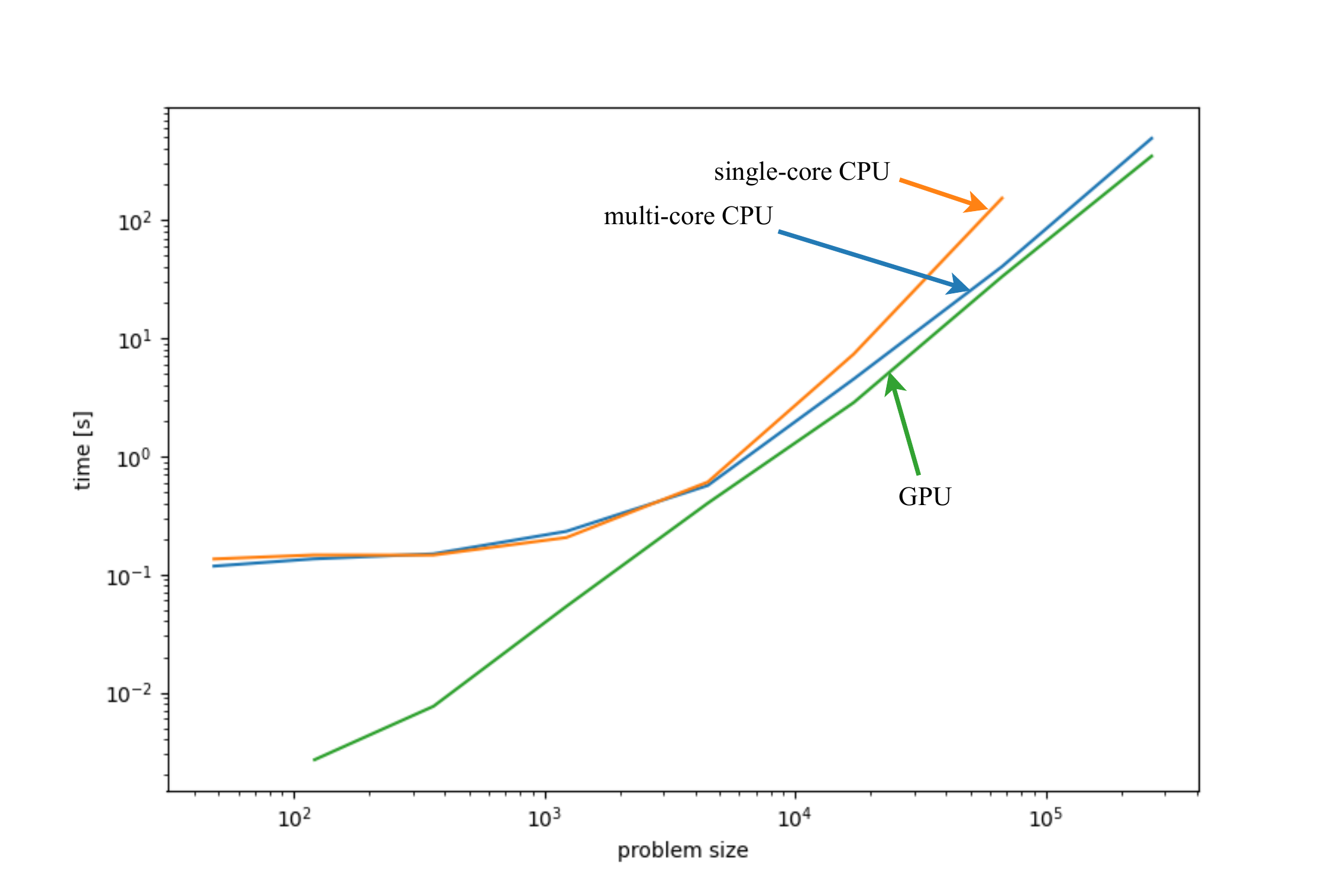}
    \caption{Comparison of CPU and GPU factorization time of sparse matrices for cubic basis functions.}
    \label{fig:res_03}
\end{figure}

All tests were performed using sparse matrices obtained from Finite Element Method (FEM) computations for the benchmark problem discussed in Section \ref{sec:formulation}.
The matrices were generated using external software and subsequently tested with two different implementations of sparse matrix solvers.
The CPU implementation was developed in Fortran, with parallelization achieved using OpenMP.
On the other hand, the GPU version was implemented in C++ and utilized CUDA for parallelization.
Performance measurements were obtained for various mesh sizes and polynomial degrees.
The reported time results are for $p = 1$ (Figure \ref{fig:res_01}), $p = 2$ (Figure \ref{fig:res_02}), and $p = 3$ (Figure \ref{fig:res_03}).
Our experiments were limited to $3 \times 10^5$ degrees of freedom for linear, quadratic, and cubic basis functions.



The comparison shown in Figures \ref{fig:res_01}-\ref{fig:res_03} highlights the substantial speed advantage of sparse matrix factorization on a shared memory graphics card, surpassing integration on 4 CPU cores by one to two orders of magnitude.

Nevertheless, it is important to note that the obtained results may have been affected by an issue related to vRAM on the GPU during our computations on Google Colab. The GPU utilized more RAM than expected, reaching its capacity for the largest matrices. We suspect this may be a Google Colab-specific problem.

Our implementation has been optimized for optimal performance, using our knowledge and expertise in optimizing summations and reductions.
However, it is worth mentioning that our current implementation does not support multi-GPU or multi-node configurations. Despite this limitation, we can still estimate scalability in such scenarios.

Adding more computational nodes may lead to efficiency drops due to communication costs, particularly in NUMA architectures \cite{numa1, numa2}. The extent of these efficiency drops can vary significantly depending on the specific multi-node, multi-GPU architecture employed.

\section{Conclusions}
\label{sec:conclusions}


We have developed a trace theory-based parallel multi-frontal direct solver, which provides efficient solutions to linear systems of equations.
The solver algorithm is expressed as a sequence of atomic tasks, each with specific dependencies between them.
We have identified the task alphabet and dependency relations for the sequential trace applied to the Finite Element Method.

Our methodology enables near-optimal scheduling on modern parallel architectures, leveraging the Foata Normal Form.
The solver execution involves multiple Foata classes, where independent tasks are computed concurrently within each class.
Synchronization barriers are employed between subsequent classes to ensure correct dependencies.
We tested an NVIDIA Tesla GPU and a multicore CPU to evaluate the solver's performance.
The methodology was validated using a simple benchmark problem for the Finite Element Method, but its applicability extends to any sparse matrix.

Looking ahead, future work can explore the application of our methodology to create state-of-the-art implementations for cutting-edge machines, such as the recently introduced NVIDIA Grace Hopper architecture \cite{grace-hopper}.
By leveraging the trace theory-based approach, we can further enhance the solver's efficiency and scalability, enabling it to tackle complex problems on advanced computing platforms.

\vspace{0.5cm}
\noindent{\bf{Acknowledgments}}
This project has received funding from The European Union’s Horizon 2020 Research and Innovation Program of the Marie Sk{\l}odowska-Curie grant agreement No. 777778, MATHROCKs.
The work of SR was partially supported by the Chilean grant ANID FONDECYT No. 3210009.

\bibliographystyle{elsarticle-num}
\bibliography{bibliography}

\end{document}